# A DATABASE FOR FIELD EXTENSIONS OF THE RATIONALS

JÜRGEN KLÜNERS AND GUNTER MALLE

ABSTRACT. We announce the creation of a database for number fields. We describe the contents and the methods of access, indicate the origin of the polynomials and formulate the aims of this collection of fields.

## 1. INTRODUCTION

We report on a database of field extensions of the rationals, its properties and the methods used to compute it. At the moment the database encompasses roughly 100,000 polynomials generating distinct number fields over the rationals, of degrees up to 15. It contains polynomials for all transitive permutation groups up to that degree, and even for most of the possible combinations of signature and Galois group in that range. Moreover, whenever these are known, the fields of minimal discriminant with given group and signature have been included. The database can be downloaded from

> www.iwr.uni-heidelberg.de/iwr/compalg/minimum/minimum.html

or from

> www.mathematik.uni-kassel.de/~malle/minimum/minimum.html

and accessed via the computer algebra system Kant [10].

One of the aims of the compilation of this database was to test the limitations of current methods for the realization of groups as Galois groups. It turned out that these methods have limitations if the signature of the resulting Galois extension is also prescribed.

## 2. GALOIS REALIZATIONS WITH PRESCRIBED SIGNATURE

Let $K/\mathbb{Q}$ be a number field of degree $n$. We denote by $r_1$ the number of real embeddings of $K$ and by $r_2$ the number of pairs of complex embeddings. Then we have $n = r_1 + 2r_2$. The pair $(r_1, r_2)$ is called the signature of $K$. The extension $K/\mathbb{Q}$ is called *totally real* if $r_2 = 0$. The solution of embedding problems often requires the knowledge of field extensions with a prescribed signature. This is one reason for the attempt to realize all groups in all possible signatures.

Now let $G$ be the Galois group of the Galois closure of $K/\mathbb{Q}$. Then for any embedding of $K$ into $\mathbb{C}$ complex conjugation is an element of $G$, that is, $G$ in its permutation representation on the conjugates of the fixed group of $K$ contains an involution with $r_1$ fixed points. Clearly this restricts the signatures which may occur for a given Galois group. It leads to the following question:

*Given a finite permutation group $G$ and a conjugacy class $C$ of involutions in $G$, does there exist a number field $K/\mathbb{Q}$ whose Galois closure has group $G$ such that the image of complex conjugation lies in class $C$?*

Obviously a positive solution to this problem would in particular solve the inverse problem of Galois theory. Serre [33] has remarked that the converse is true at least for totally real extensions:

**Proposition 1.** *(Serre) If all finite groups occur as Galois groups over $\mathbb{Q}$, then all finite groups occur as Galois groups of totally real extensions of $\mathbb{Q}$.*

*Proof.* Let $G$ be a finite group. We use a special case of a result of Haran and Jarden [20], Cor. 6.2: There exists a finite group $\tilde{G}$ and an epimorphism $\phi : \tilde{G} \to G$ having all involutions of $\tilde{G}$ in the kernel. Indeed, assume that $G$ is generated by $r$ elements, and let $\psi : F_r \to G$ be a corresponding epimorphism from the free profinite group $F_r$ of rank $r$ onto $G$. Then $S := F_r \setminus \ker(\psi)$ is a compact subset which is invariant under conjugation. Hence $S^2 := \{g^2 \mid g \in S\}$ is compact and invariant under conjugation as well and $1 \notin S^2$. Thus there exists a normal subgroup $\tilde{N} \le F_r$ of finite index with $\tilde{N} \cap S^2 = \emptyset$. Defining $N := \tilde{N} \cap \ker(\psi)$ we have $N \cap S^2 = \emptyset$ and $N$ has finite index in $F_r$. Then $\tilde{G} := F_r/N$ with $\phi : \tilde{G} \to G$ induced by $\psi$ is as required.

Now let $K/\mathbb{Q}$ be a Galois extension with group $\tilde{G}$ and $\phi : \tilde{G} \to G$ as above with $H = \ker(\phi)$. Assume that $K^H/\mathbb{Q}$ is not totally real. Then some involution of $G$ acts as complex conjugation. But by construction of $\tilde{G}$ this involution lifts to an element of order bigger than 2 in $\tilde{G} = \mathrm{Gal}(K/\mathbb{Q})$, contradicting the fact that complex conjugation has order 2. Thus $K^H/\mathbb{Q}$ is a totally real realization for $G$. □

For the explicit construction of fields with given signature, we may distinguish two cases. In the solvable case, class field theory may be used as in the general inverse problem. The construction of extensions with non-solvable groups usually is done via the rigidity method. But this seems less adapted to the case where in addition the signature is prescribed. In fact, Serre [32], p.91, has shown, that rigidity with three branch points never gives totally real Galois extensions for groups $G \ne \mathfrak{S}_3$. At the moment we are reduced to ad hoc methods for constructing extensions with arbitrary signature.

2.1. **Symmetric groups.** Let's first treat the symmetric groups. We propose an even stronger statement.

**Proposition 2.** *Let $n \in \mathbb{N}$, $0 \le k \le n/2$, and $f_i \in \mathbb{Q}_{p_i}[X]$ ($i = 1, \ldots, r$) polynomials of degree $n$, where $p_i \ne p_j$ for $i \ne j$. Then there exist infinitely many number fields $K/\mathbb{Q}$ with Galois group $\mathfrak{S}_n$, signature $(n - 2k, k)$ and such that $K \otimes \mathbb{Q}_{p_i} \cong \mathbb{Q}_{p_i}[X]/(f_i)$ for $i = 1, \ldots, r$.*

*Proof.* Let $g_0 \in \mathbb{Z}[X]$ be a separable polynomial with $n - 2k$ real and $k$ pairs of complex zeros, for example the polynomial $\prod_{i=1}^{n-2k}(X - i) \prod_{i=1}^{k}((X - i)^2 + 1)$. By the main theorem on elementary symmetric functions and Hilbert's irreducibility theorem there exist irreducible polynomials with group $\mathfrak{S}_n$ arbitrarily close to $g_0$ (for example with respect to the metric induced by taking the maximal absolute value of the coefficients). To find such a polynomial constructively, choose three further primes $p_{r+1}, \ldots, p_s$, $s := r + 3$. Furthermore, let $g_i \in \mathbb{Z}[X]$, $i = 1, \ldots, r$, be polynomials such that $\mathbb{Q}_{p_i}[X]/(f_i) \cong \mathbb{Q}_{p_i}[X]/(g_i)$, $i = 1, \ldots, r$, and $g_{r+1}, \ldots, g_s \in \mathbb{Z}[X]$ separable such that the only non-linear irreducible factor of the reduction $g_i$ (mod $p_i$) ($r + 1 \le i \le s$) has degree $n$, $n - 1$, 2 respectively.

Write
$$g_i = \sum_{j=0}^{n} a_{i,j} X^j \text{ for } i = 0, \ldots, s.$$

By the weak approximation theorem [25, Theorem 1.11] we may choose $b_0, \ldots, b_n \in \mathbb{Q}$ such that
$$f := \sum_{i=0}^{n} b_i X^i \in \mathbb{Q}[X]$$
satisfies the following two conditions:
- $|b_j - a_{0,j}|$ $(j = 0, \ldots, n)$ is sufficiently small so that $f$ has the same signature as $g_0$,
- $|b_j - a_{i,j}|_{p_i}$ are sufficiently small so that $\mathbb{Q}[X]/(f) \otimes \mathbb{Q}_{p_i} \cong \mathbb{Q}_{p_i}[X]/(g_i)$ for $i = 1, \ldots, s$. (Lemma of Krasner, [25, Proposition 5.5]).

Then $f$ has the same signature as $g_0$, and $\mathbb{Q}[X]/(f) \otimes \mathbb{Q}_{p_i} \cong \mathbb{Q}_{p_i}[X]/(f_i)$, $i = 1, \ldots, r$, are as required. Finally, factorization modulo $p_{r+1}, \ldots, p_s$ shows that $\text{Gal}(f)$ is 2-fold transitive, hence primitive, and contains a transposition. By a theorem of Jordan this implies that $\text{Gal}(f) = \mathfrak{S}_n$.
By varying the additional primes, respectively by enlarging the set of primes, we may clearly obtain infinitely many examples. □

2.2. **Alternating groups.** The case of alternating groups is already more complicated, and it is the only other case which we can solve uniformly. We first rephrase a result of Mestre into a universal lifting property (see [2] for definitions and other results in this direction):

**Theorem 3.** *The group $\mathfrak{A}_n$ has the universal lifting property over fields of characteristic 0. More precisely, if $K$ is a field of characteristic 0 and $g(X) \in K[X]$ a separable polynomial with square discriminant then there exists a polynomial $f(t, X) \in K(t)[X]$ generating a regular Galois extension of $K(t)$ with group $\mathfrak{A}_n$, $n \geq 3$, such that the splitting fields of $g(X)$ and $f(0, X)$ coincide.*

*Proof.* Let $g(X) \in K[X]$ of degree $n \geq 3$ with square discriminant. First assume that $n$ is odd. Then by the result of Mestre (see [24], IV.5.12) there exists a polynomial $h(X) \in K(X)$ of degree $n - 1$ such that $f(t, X) := g(X) - th(X) \in K(t)[X]$ has Galois group $\mathfrak{A}_n$ over $K(t)$.
Now assume that $n$ is even. Replacing $X$ by $X - a$ for a suitable $a \in K$ we may assume that $g_1(X) := Xg(X)$ is separable. Since $g_1$ has again a square discriminant and odd degree, by the first part there exists a polynomial $f_1(t, X) = g_1(X) - th(X)$ with group $\mathfrak{A}_{n+1}$. Note that this implies that $h(0) \neq 0$. By [24], IV.5.12(b), the polynomial
$$\tilde{f}(t, X) := (g_1(X)h(t) - g_1(t)h(X))/(X - t) \in K(t)[X]$$
has group $\mathfrak{A}_n$. Moreover, $\tilde{f}(0, X) = Xg(X)h(0)/X = h(0)g(X)$ is a non-zero scalar multiple of $g(X)$, so $f(t, X) := \tilde{f}(t, X)/h(0)$ has all the required properties. □

Note that the signature of a field with even Galois group is necessarily of the form $(n - 2k, k)$ with $k$ even. For alternating groups, all these signatures can be realized over $\mathbb{Q}$.

**Corollary 4.** *Let $n \in \mathbb{N}$ and $0 \leq k \leq n/2$ even. Then there exist infinitely many number fields $K/\mathbb{Q}$ with Galois group $\mathfrak{A}_n$ and signature $(n - 2k, k)$.*

*Proof.* For $1 \leq i \leq k/2$ let $u_i(X) \in \mathbb{Z}[X]$ be distinct totally complex polynomials of degree 4 with Galois group $\mathfrak{A}_4$, for example $u_i(X) = (X-i)^4 - 7(X-i)^2 - 3(X-i) + 1$. Then

$$g(X) := \prod_{i=1}^{n-2k} (X-i) \prod_{i=1}^{k/2} u_i(X)$$

is separable with square discriminant and signature $(n-2k, k)$. By Theorem 3 there exists a polynomial $f(t, X) \in \mathbb{Q}(t)[X]$ with Galois group $\mathfrak{A}_n$ such that $f(0, X) = g(X)$. Since $g$ is separable, for any $t_0$ close to 0 the specialization $f(t_0, X)$ has the same signature as $g$. By the Hilbert irreducibility theorem there exist infinitely many such $t_0$ for which the Galois group is preserved under specialization. □

Note that for symmetric and alternating groups the conjugacy classes of involutions are parameterized by the cycle types, so the preceding results show that any involution in an alternating or symmetric group can occur as complex conjugation in a Galois extension of the rationals.

2.3. **Further simple groups.** The non-abelian simple groups with faithful permutation representations of degree at most 15 are $L_2(7), L_2(8), L_2(11), M_{11}, M_{12}, L_3(3), L_2(13)$ and the alternating groups. For the groups $L_2(7), L_2(11), M_{11}$ and $M_{12}$ totally real realizations were found in [23], by constructing the Hurwitz spaces for certain $n$-tuples of conjugacy classes, where $n \geq 4$. These constructions involve a considerable amount of calculation and seem to be restricted to small degree. At the moment we are not aware of any totally real extensions of $\mathbb{Q}$ with group $L_2(8), L_3(3)$ or $L_2(13)$, nor with the almost simple groups $P\Gamma L_2(8)$, $PGL_2(11)$ or $PGL_2(13)$.

## 3. How to construct polynomials

In this section we give a short overview about the methods we used to construct the polynomials contained in the database.

3.1. **Methods from the geometry of numbers.** Let $K$ be a number field of degree $n$ with absolute discriminant $D$. For $\alpha \in K$ we denote by $\alpha = \alpha_1, \ldots, \alpha_n \in \mathbb{C}$ the conjugates of $\alpha$ and define $T_2(\alpha) := \sum_{i=1}^n |\alpha_i^2|$. Now Theorem 6.4.2 in [7] (attributed to Hunter) states that there exists an algebraic integer $\alpha \in K \setminus \mathbb{Q}$ such that $T_2(\alpha) \leq B$, where $B$ only depends on $n, D$. This can be used to derive bounds for the coefficients of the characteristic polynomial of a primitive element of $K$. A description of this method can be found in [8, Section 9.3]. In the case that all conjugates $\alpha_1, \ldots, \alpha_n$ are real we have used a slightly different approach.

Let $f(X) \in \mathbb{Z}[X]$ be a totally real separable polynomial of degree $n$ (that is, the stem field of $f$ has $n$ different real embeddings). Then all derivatives of $f$ are also totally real and separable. Conversely given a totally real polynomial $g(X) \in \mathbb{Z}[X]$ of degree $n-1 \geq 2$ there are only finitely many totally real polynomials $f(X) \in \mathbb{Z}[X]$ such that $f' = g$. Moreover the constant terms of such polynomials $f$ consist of all integers in an interval $I$ which can be computed from $g$. Indeed, denote by $\alpha_1 < \ldots < \alpha_{n-1}$ the (different real) roots of $g$ and let $f_0$ denote any integral polynomial with derivative $f_0' = g$. Assume for definiteness that the highest coefficient of $g$ is positive. Denote by

$$m := \max\{f_0(\alpha_{n-1-2i}) \mid 0 \leq i \leq (n-2)/2\}$$

the maximum of the minima of $f_0$, and by

$$M := \min\{f_0(\alpha_{n-2-2i}) \mid 0 \leq i \leq (n-3)/2\}$$

the minimum of the maxima. Then clearly $f_0 - c$ is totally real if and only if $c \in I := \{\alpha \in \mathbb{R} \mid m \leq \alpha \leq M\}$.

The above considerations show the following (see also [8, p. 448], for example):

**Lemma 5.** *For fixed $a_n, a_{n-1}, a_{n-2} \in \mathbb{Z}$ there exist only finitely many totally real polynomials $f(X) = \sum_{i=0}^{n} a_i X^i \in \mathbb{Z}[X]$, and these may be enumerated effectively.*

Indeed, such polynomials can only exist if $f^{(n-2)}$ is totally real. Since $f^{(n-2)}$ of degree 2 is completely determined by $a_n, a_{n-1}, a_{n-2}$, there are only finitely many possibilities for $f^{(n-3)}$, and now induction proves the assertion.

By the theorem of Hunter any primitive extension $K$ of $\mathbb{Q}$ of degree $n$ can be generated by a monic polynomial $f(X) = \sum_{i=0}^{n} a_i X^i \in \mathbb{Z}[X]$ of degree $n$ such that $0 \leq a_{n-1} \leq n/2$ and with $T_2$-norm bounded by a function in the discriminant $d(K)$. Moreover the $T_2$-norm bounds the third highest coefficient $a_{n-2}$.

Hence Lemma 5 can be used to enumerate totally real fields of bounded discriminants. It seems that this strategy produces much fewer polynomials to be considered, as compared to the approach which first tries to bound the discriminant and then to sieve for totally real polynomials. For example, in the case of totally real degree 8 extensions (see Theorem 13), only 869062 polynomials were produced and had to be processed further. (Among the corresponding fields only 4896 had Galois group different from $\mathfrak{S}_8$.)

3.2. **Specializing from polynomials over $\mathbb{Q}(t)$.** Let $G$ be a finite group. We call the field extension $K/\mathbb{Q}(t)$ a $G$-realization, if it is Galois with group $G$ and regular, which means that $\mathbb{Q}$ is algebraically closed in $K$. When a group has a $G$-realization over $\mathbb{Q}$, it is an immediate consequence that there exist infinitely many disjoint number fields $L/\mathbb{Q}$ with Galois group $G$: Suppose we have a polynomial $f \in \mathbb{Q}(t)[X]$ such that the splitting field of $f$ is a regular extension with Galois group $G$. By specializing $t$ to $a \in \mathbb{Q}$ we get that $\text{Gal}(f(a, X))$ is a subgroup of $G$. Hilbert's irreducibility theorem states that $\text{Gal}(f(a, X)) = G$ for infinitely many $a \in \mathbb{Q}$. See for example [32, Section 4.6] for a method to find infinitely many $a \in \mathbb{Q}$ with that property. This allows to construct polynomials with Galois group $G$ over $\mathbb{Q}$ when we have an explicit polynomial $f \in \mathbb{Q}(t, X)$. In some lucky cases we are able to get proper subgroups of $G$.

3.3. **Methods from class field theory.** Suppose we want to construct a polynomial $f$ such that $\text{Gal}(f) = G$ for some permutation group $G$. Furthermore suppose that in a corresponding field extension the stem field $N$ of $f$ has a subfield $L$ such that $N/L$ is an Abelian extension with Galois group $A$. Then we can try the following approach. The Galois group of (the splitting field of) $L$ can be determined group theoretically and is denoted by $H$. Given a field $L$ with Galois group $H$ we generate relative Abelian extensions with Galois group $A$ using class field theory. The Galois groups over $\mathbb{Q}$ of such extensions are subgroups of the wreath product $A \wr H$. Experiments show that most of the computed fields have the wreath product or the direct product as Galois group. But we also get other Galois groups. One advantage of this method is that we are able to control the field discriminants of the computed fields. Therefore we can prove minimal discriminants for such groups. E.g. this was applied successfully to degree 8 fields having a degree 4 subfield [9].

For a complete description we refer the reader to [8, Section 9.2]. We remark that we have used the class field algorithm described in [15] and implemented in [10]. Cohen [8, Theorem 9.2.6] remarks that the class field methods can be extended to fields where the Galois group of $N/L$ is a dihedral group of order $2n$, where $n$ is odd. C. Fieker and the first author [16] can extend this method to the case where $N/L$ is a Frobenius group with Abelian kernel. For example this applies to the Frobenius groups $Z_l \rtimes Z_p$, where $p$ is prime and $p \mid l-1$.

3.4. **Embedding obstructions.** Suppose we want to construct a field extension of degree 4 with cyclic group $Z_4$ applying the methods of the preceding paragraph and take $L := \mathbb{Q}(\sqrt{-1})$. Then we will find out that there are no extensions $N/L$ such that $\mathrm{Gal}(N/\mathbb{Q}) \cong Z_4$. It would be nice to know in advance whether $L$ is a good choice or not. Let $K$ be a number field and $L/K$ be a finite field extension with Galois group $H$ and

$$1 \to U \to G \to H \to 1$$

be an exact sequence of groups. Then a field $N/L$ is called a proper solution of the embedding problem, if $\mathrm{Gal}(N/K) \cong G$. For the general theory we refer the reader to [24, Chapter IV]. Here we restrict ourselves to the special case with kernel $U \cong Z_2$. Then $U$ is a subgroup of the center of $G$ and we have the following result [24, IV.7.2].

**Proposition 6.** *Let $N = L(\sqrt{\alpha})$ with $\alpha \in L$ be a proper solution of the given embedding problem with kernel $Z_2$. Then all solution fields are of the form $N_a := L(\sqrt{a\alpha})$ with $a \in K^\times$.*

Furthermore we get a local-global principle. Let $L/K$ be a number field with Galois group $H$ and suppose we have the embedding problem

$$1 \to Z_2 \to G \to H \to 1.$$

Denote by $\mathbb{P}_K$ the set of prime ideals of $\mathcal{O}_K$ including the infinite ones. For $\mathfrak{p} \in \mathbb{P}_K$ and $\mathfrak{P}$ a prime ideal of $\mathcal{O}_L$ lying over $\mathfrak{p}$ we denote by $L_\mathfrak{P}/K_\mathfrak{p}$ the corresponding local extension. We write $\bar{H}$ for the Galois group of $L_\mathfrak{P}/K_\mathfrak{p}$. We get the following induced embedding problem:

$$1 \to Z_2 \to \bar{G} \to \bar{H} \to 1.$$

This embedding problem has a solution if it has a proper solution or if the exact sequence is split (see [24, p. 265] for the general definition of "solution").

**Proposition 7.** *Let $L/K$ be a finite extension with Galois group $H$. Then the embedding problem $1 \to Z_2 \to G \to H \to 1$ has a proper solution if and only if the induced embedding problems have a solution for all $\mathfrak{p} \in \mathbb{P}_K$ with one possible exception.*

*Proof.* The theorem follows from [24, Cor. IV.10.2] and the subsequent remark and a theorem of Ikeda [24, Th. IV.1.8]. Recall that split embedding problems with Abelian kernel have proper solutions [24, Th. IV.2.4]. □

In our special case with kernel $Z_2$ it is easy to see that the induced embedding problems have solutions for all $\mathfrak{p}$ which are unramified in $L$ or which have odd ramification index in $L$. If an infinite prime $\mathfrak{p}$ is ramified the induced embedding problem is solvable if and only if it is split.

These results give us a practical method to check if an embedding problem with kernel $Z_2$ has a proper solution. If there exists a solution field of the embedding problem it remains to compute such a field.

**Proposition 8.** *Let $N = L(\sqrt{\alpha})$ be a proper solution of an embedding problem with kernel $Z_2$. Let $S \subset \mathbb{P}_K$ be a finite subset containing all prime ideals with even ramification index in $L/K$, all infinite primes, and all prime ideals lying above $2\mathbb{Z}$. Furthermore assume that $S$ contains enough prime ideals to generate the class group of $\mathcal{O}_K$. Then there exists a proper solution $\tilde{N}/L$ which is unramified outside $\tilde{S}$, where $\tilde{S} := \{\mathfrak{P} \in \mathbb{P}_L \mid \mathfrak{P} \supseteq \mathfrak{p} \text{ for some } \mathfrak{p} \in S\}$.*

*Proof.* Denote by $\hat{S}$ the set of all prime ideals in $\mathcal{O}_L$ which are ramified in $N$ and are not contained in $\tilde{S}$. All prime ideals in $\hat{S}$ are tamely ramified. Furthermore, if $\mathfrak{P} \in \hat{S}$ it follows that all conjugate prime ideals are contained in $\hat{S}$ as well. Define $\mathfrak{a}$ to be the product of all prime ideals contained in $\hat{S}$. We get that $\mathfrak{a} = \mathfrak{b}\mathcal{O}_L$, where $\mathfrak{b}$ is a square-free ideal in $\mathcal{O}_K$. Then there exist $\mathfrak{p}_1, \ldots, \mathfrak{p}_r \in S$ and $e_1, \ldots, e_r \in \mathbb{N}$ such that $\mathfrak{b}\mathfrak{p}_1^{e_1} \cdots \mathfrak{p}_r^{e_r}$ is a principal ideal in $\mathcal{O}_K$ with generator $b$. Then $N_b := L(\sqrt{b\alpha})$ is a proper solution unramified outside $\tilde{S}$. $\square$

Since there are only finitely many relative quadratic extensions unramified outside a finite set, the above furnishes a method to explicitly compute a solution. We remark that in the case $K = \mathbb{Q}$ the condition about the infinite primes can be dropped. In the case that $L$ is totally real and $L(\sqrt{\alpha})$ is totally complex (both extensions are normal over $\mathbb{Q}$) the field $L(\sqrt{-\alpha})$ is totally real.

Let us give a few examples how the solvability in the $p$-adic case can be decided.

**Example 1.**   1. *A degree 2 extension $L/\mathbb{Q}$ is embeddable into a $Z_4$ extension if and only if $L$ is totally real and all odd primes $p$ which are ramified in $L$ are congruent $1 \bmod 4$.*
  2. *Let $L/\mathbb{Q}$ be an extension with Galois group $V_4$. Then $L$ is embeddable into a $Q_8$ extension if and only if $L$ is totally real and all odd primes $p$ which are ramified in $L$ have the property that $p \equiv 1 \bmod 4$ if and only if $p$ has odd inertia degree in $L$.*
  3. *Let $L/\mathbb{Q}$ be an extension with Galois group $\mathrm{L}_2(l)$, where $l$ is a prime with $l \equiv 3 \bmod 8$ or $l \equiv 5 \bmod 8$. Then $L$ is embeddable into an $\mathrm{SL}_2(p)$ extension if and only if $L$ is totally real and all odd primes $p$ which are ramified in $L$ have the property that $p \equiv 1 \bmod 4$ if and only if $p$ has odd inertia degree in $L$ [5].*

The following example is more complicated and demonstrates most of the effects which may occur.

There exists a subdirect product $G = \mathrm{SL}_2(3) \times_{\mathfrak{A}_4} [4^2]3$ with a faithful transitive permutation representation of degree 12, usually denoted $12T_{57}$. As we have noted in [21, 4.1] in order to construct an extension with this group we have to find an $\mathfrak{A}_4$-extension which is embeddable both into an $\mathrm{SL}_2(3)$-extension and into a $[4^2]3$-extension. For $p \neq 2$ the possible non-trivial local Galois groups of an $\mathfrak{A}_4$-extension are $Z_2, Z_3, Z_2 \times Z_2$. Let $E/\mathbb{Q}_p$, $p \neq 2$, be a $p$-adic field. If the local Galois group is totally ramified with Galois group $Z_2$ we get that both local embedding problems are solvable if $p \equiv 1 \bmod 4$. If the local Galois group is $Z_2 \times Z_2$ it cannot be a totally ramified extension ($p \neq 2$, Abhyankar's lemma [25, p. 236]). In this case the embedding problem into $\mathrm{SL}_2(3)$ can only be solved when $p \equiv 3 \bmod 4$. But then the 4th roots of unity are not contained in $\mathbb{Q}_p$ and the embedding problem into $[4^2]3$ cannot be solved. Therefore we get: Let $L/\mathbb{Q}$ be an extension which is embeddable into a $12T_{57}$ extension. Then $L$ is totally real and odd primes $p$ which are ramified have inertia degree 1 and satisfy $p \equiv 1 \mod 4$. The converse is true when $L$ is unramified in 2 or the degree of the completion at 2 has degree divisible by 3.

**Proposition 9.** *Let $L/\mathbb{Q}$ be an extension with Galois group $\mathfrak{A}_4$. Then $L$ is embeddable into a $12T_{57}$ extension if and only if the following holds:*

1. *$L$ is totally real.*
2. *If $p \neq 2$ is a ramified prime in $L$ then $p \equiv 1 \mod 4$ and $p$ has inertia degree 1 in $L$.*
3. *If 2 is ramified, then the corresponding embedding problem for $p = 2$ is solvable.*

*Denote by $M$ the subfield of $L$ which has Galois group $Z_3$. Suppose that $L$ is embeddable into a $12T_{57}$ extension. Denote by $S$ the set of prime ideals in $\mathcal{O}_M$ containing all prime ideals above $2\mathbb{Z}$, all infinite primes, all prime ideals which are ramified in $L$, and enough prime ideals to generate the class group of $\mathcal{O}_M$. Then there exists a $12T_{57}$ extension containing $L$ which is unramified outside $\tilde{S}$, where $\tilde{S} := \{\mathfrak{P} \in \mathbb{P}_L \mid \mathfrak{P} \supseteq \mathfrak{p} \text{ for some } \mathfrak{p} \in S\}$.*

*Proof.* The first part of the theorem is already proved. We can solve the corresponding embedding problems independently. For the $\mathrm{SL}_2(3)$-part we can apply Proposition 8. Denote by $K$ one of the degree 6 subfields of $L$. As noted in [21, 4.1] the embedding problem into $[4^2]3$ is solvable if and only if $K/M$ is embeddable into a $Z_4$-extension. Therefore we can again apply Proposition 8. □

If 2 is ramified we cannot decide the solvability of the embedding problem just by looking at the ramification behaviour. We have to determine if $K/M$ is embeddable into a $Z_4$ extension, which is the case if and only if $-1$ is a norm in $K/M$. This can be decided by applying the methods described in [1].

3.5. **Computing polynomials from other representations.** Suppose we want to compute polynomials for a permutation group which already has a faithful representation on fewer points, that is, we want to construct a different stem field of a given Galois extension. In [21, 3.3] we have described how to compute such polynomials when we know a polynomial belonging to the other representation. In this paper we strive to control the discriminant of these fields. The proof of the following theorem can be found in [22, Proposition 6.3.1].

**Theorem 10.** *Let $N/K$ be a normal extension with Galois group $G$ and $L$ be the fixed field of a subgroup $H$ of $G$. Let $\mathfrak{P} \neq (0)$ be a prime ideal of $\mathcal{O}_N$ with ramification index $e$ and $\mathfrak{p} := \mathfrak{P} \cap \mathcal{O}_K$. Denote by $D_{\mathfrak{P}}$ and $I_{\mathfrak{P}}$ the decomposition group and inertia group, respectively. Let $R_H := \{g_1, \ldots, g_m\}$ be a system of representatives of the double cosets of $H$ and $D_{\mathfrak{P}}$ in $G$, i.e. $G = \dot{\bigcup}_{i=1}^{m} H g_i D_{\mathfrak{P}}$. Then*

1. *The prime divisors of $\mathfrak{p}$ in $\mathcal{O}_L$ are $\mathfrak{p}_i := g_i \mathfrak{P} \cap \mathcal{O}_L$ for $1 \leq i \leq m$.*
2. *$\mathfrak{p}\mathcal{O}_L = \prod_{i=1}^{m} \mathfrak{p}_i^{e_i}$, where $e_i := \frac{e}{|g_i I_{\mathfrak{P}} g_i^{-1} \cap H|}$.*

**Corollary 11.** *Suppose that $\mathfrak{P}$ is not wildly ramified over $\mathfrak{p}$. In this case we denote by $\pi$ a generator of the cyclic group $I_{\mathfrak{P}}$. Then $v_{\mathfrak{p}}(\mathrm{disc}(L/K)) = \mathrm{ind}(\pi)$, where $\mathrm{ind}(\pi) := [G:H]-$ number of orbits of $\pi$ on $G/H$.*

*Proof.* Suppose that $\mathfrak{p}\mathcal{O}_L = \prod_{i=1}^{m} \mathfrak{p}_i^{e_i}$. In the case of tamely ramified extensions we get that $v_{\mathfrak{p}}(\mathrm{disc}(L/K)) = \sum_{i=1}^{m} f_i(e_i - 1)$, where $f_i$ denotes the degree of the residue field extension $(\mathcal{O}_L/\mathfrak{p}_i)/(\mathcal{O}_K/\mathfrak{p})$. Obviously this formula does not depend on the number of primes lying above $\mathfrak{p}$ or their inertia degrees, only on the index of $\pi$. □

**Example 2.** *Let $G = \mathrm{L}_2(7)$, the second smallest non-abelian simple group. The following table illustrates the assertion of the previous Corollary. The two columns give the cycle types of elements of $G$ in the transitive degree 7 and degree 8 representations, respectively. This allows to compare the contribution to the discriminant by tamely ramified prime ideals.*

| $n=7$ | $n=8$ |
|---|---|
| $1^7$ | $1^8$ |
| $1^3 \cdot 2^2$ | $2^4$ |
| $1 \cdot 3^2$ | $1^2 \cdot 3^2$ |
| $1 \cdot 2 \cdot 4$ | $4^2$ |
| $7$ | $1 \cdot 7$ |

*We can see that independently of the cycle type the discriminant in the degree 8 representation remains at least the same as in the degree 7 representation, in case of tame ramification. A case by case study shows that the same is true when wild ramification occurs. This opens a way to determining the smallest fields of degree 8 with Galois group $\mathrm{L}_2(7)$ by computing enough fields of degree 7 with the corresponding Galois group.*

## 4. Minimal discriminants

**4.1. Results known to date.** One goal of our database is to provide fields with small (absolute value of the) discriminant for each Galois group and signature. In small degrees it is even possible to determine the field(s) with smallest discriminant. Let's comment on the present state of knowledge in this area (which is restricted to degrees less than 10).

It is very easy to enumerate the discriminants of quadratic fields. Belabas [3] gives a very efficient algorithm to enumerate cubic number fields. For higher degrees methods from the geometry of numbers and class field theory are applied.

In [6] all quartic fields with absolute discriminant smaller than $10^6$ are enumerated. There are huge tables of the smallest quintic fields due to [31]. These tables are sufficient to extract the smallest discriminants for all Galois groups and classes of involutions for degrees 4 and 5.

The general enumeration methods are not powerful enough to give the minima for all Galois groups in degree 6. The minimal discriminants for all signatures of degree 6 are computed in [30]. [26, 17, 18, 19] have finished the computation of minimal discriminants of all signatures and all primitive Galois groups of degree 6. [27, 4] compute the minimal fields for imprimitive groups of degree 6. This yields enough information to determine the minimal fields for all groups and all conjugacy classes of that degree.

In degree 7 the minimal fields of each signature are known due to [12, 14, 29]. This covers all signatures of the symmetric groups. We complete the determination in degree 7 by proving the following:

**Theorem 12.** *The minimal discriminants for the possible pairs $(G, r_1)$ of Galois group $G$ and number of real places $r_1$ in degree 7 are as shown in the following table.*

Minimal discriminants in degree 7

| $G$ | $r_1 =$ | 1 | 3 | 5 | 7 |
|---|---|---|---|---|---|
| 7 | | — | — | — | 594823321 |
| 7 : 2 | | −357911 | — | — | 192100033 |
| 7 : 3 | | — | — | — | 1817487424 |
| 7 : 6 | | −38014691 | — | — | 12431698517 |
| $L_3(2)$ | | — | 2007889 | — | 670188544 |
| $\mathfrak{A}_7$ | | — | 3884841 | — | 988410721 |
| $\mathfrak{S}_7$ | | −184607 | 612233 | −2306599 | 20134393 |

*Proof.* The fields generated in [12, 14, 29] are sufficient to prove the minimal discriminants for all signatures of the symmetric group and the non totally real dihedral case. The minima for $\mathfrak{A}_7$ and $L_3(2)$ are found by using methods from the geometry of numbers (see Section 3.1) using the fact that the discriminant has to be a square. The minimal discriminant for the cyclic case can easily be determined using the theorem of Kronecker-Weber and the fact that a ramified prime $p$ must either be 7 or $p \equiv 1 \mod 7$. All the other groups are Frobenius groups, where we can apply class field theory as described in [16] to prove the minima. □

The polynomial
$$X^7 - 2X^6 - 7X^5 + 11X^4 + 16X^3 - 14X^2 - 11X + 2$$
generates a totally real $\mathfrak{A}_7$-extension with minimal discriminant, while
$$X^7 - 8X^5 - 2X^4 + 15X^3 + 4X^2 - 6X - 2$$
generates one of the two totally real $L_3(2)$-extension with minimal discriminant (the other one is arithmetically equivalent to the first one, which means that these two non-isomorphic fields have the same Dedekind $\zeta$-function).

The smallest totally real octic number field is computed in [28]. Diaz y Diaz [13] determined the smallest totally complex octic number field. To the best of our knowledge the smallest totally real octic field with symmetric Galois group was previously unknown. The following theorem can be proved using the methods of Section 3.1.

**Theorem 13.** *The minimal discriminant for a totally real primitive field of degree 8 is given by $d = 483345053$. The corresponding extension is unique up to isomorphism, with Galois group $\mathfrak{S}_8$, generated by the polynomial*
$$X^8 - X^7 - 7X^6 + 4X^5 + 15X^4 - 3X^3 - 9X^2 + 1.$$

For imprimitive octic fields with a quartic subfield [9] compute huge tables using class field theory which cover all imprimitive groups and all possible signatures such that the corresponding field has a quartic subfield. These tables are not sufficient to find all minimal fields of that shape such that complex conjugation lies in a given class of involutions. In [16] the minima for octic fields having a quadratic subfield are computed.

It remains to say something about primitive groups in degree 8. In the following table we give the primitive groups and the smallest discriminants we know. If there is no $\leq$ or $\geq$ sign this means that this entry is proven to be minimal. The totally real $\mathfrak{S}_8$ case is already proved in Theorem 13. The minima for the groups $8T_{25}$ and $8T_{36}$ are proved in [16].

Minimal discriminants of primitive groups in degree 8

| $G$ | $r_1 =$ 0 | 2 | 4 | 6 | 8 |
|---:|---:|---:|---:|---:|---:|
| $8T_{25}$ | 594823321 | — | — | — | 9745585291264 |
| $8T_{36}$ | 1817487424 | — | — | — | 6423507767296 |
| $8T_{37}$ | $\leq 37822859361$ | — | — | — | $\leq 8165659002209296$ |
| $8T_{43}$ | $\leq 418195493$ | $\geq -1997331875$ | — | — | $\leq 312349488740352$ |
| $8T_{48}$ | $\leq 32684089$ | — | $\leq 351075169$ | — | $\leq 81366421504$ |
| $8T_{49}$ | $\leq 20912329$ | — | $\leq 144889369$ | — | $\leq 46664208361$ |
| $8T_{50}$ | $\leq 1282789$ | $\geq -4296211$ | $\leq 15908237$ | $\geq -65106259$ | 483345053 |

If we knew enough fields of degree 7 with Galois group $L_2(7)$ it would be possible to compute the minima for the groups $8T_{37} \cong L_2(7)$ and $8T_{48} \cong 2^3.L_2(7)$.

Diaz y Diaz and Olivier [11] have applied a relative version of the geometry of numbers methods to compute tables of imprimitive fields of degree 9. These tables do not cover all imprimitive Galois groups of that degree.

## 5. The database

In this section we report on the content of the database. As mentioned in the introduction it contains about 100,000 polynomials generating distinct number fields over the rationals. Especially in smaller degrees (up to degree 5) there already exist much larger tables of number fields covering all fields up to a given discriminant bound. It is not very surprising that most of these fields have symmetric Galois group. The aim of our database is different. We want to cover all groups. More precisely we want to look at the following problems of increasing difficulty:

1. For each transitive group $G$ find a polynomial $f \in \mathbb{Z}[x]$ such that $\mathrm{Gal}(f) = G$.
2. For each transitive group $G$ and each class $C$ of involutions find a polynomial $f \in \mathbb{Z}[x]$ such that $\mathrm{Gal}(f) = G$ and complex conjugation lies in class $C$.
3. For each transitive group $G$ and each class $C$ of involutions find a polynomial $f \in \mathbb{Z}[x]$ such that $\mathrm{Gal}(f) = G$ and complex conjugation lies in class $C$ and the stem field $K$ of $f$ has minimal absolute discriminant subject to these restrictions.

We have a positive answer to problem 1 for all transitive groups up to degree 15, as shown in [21]. Problem 2 is already much more difficult. Let us first look at a slightly easier variant of problem 2. Here we only ask that complex conjugation covers all cycle types of involutions in $G$. The easier problem has a positive answer for all transitive groups but the following possible exceptions:

| Group | Number of real zeroes |
|---:|---:|
| $9T_{27} = L_2(8)$ | 9 |
| $9T_{30} = P\Gamma L_2(8)$ | 9 |
| $12T_{218} = PGL_2(11)$ | 12 |
| $13T_7 = L_3(3)$ | 13 |
| $13T_8 = \mathfrak{A}_{13}$ | 5,9,13 |
| $14T_{30} = L_2(13)$ | 14 |
| $14T_{39} = PGL_2(13)$ | 14 |
| $14T_{62} = \mathfrak{A}_{14}$ | 6 |
| $15T_{103} = \mathfrak{A}_{15}$ | 7,11,15 |

The missing signatures for the alternating groups are only a practical problem as we have proved in Theorem 3. In all the other cases the missing signature is the totally real one; we don't even know a theoretical argument that such an extension should exist.

Let us come back to problem 2. Write $N := N_{\mathfrak{S}_n}(G)$ for the normalizer in $\mathfrak{S}_n$ of $G \leq \mathfrak{S}_n$. Let $L/K$ be an extension of degree $n$ generated by a polynomial $f$ such that $G$ is the Galois group of the Galois closure of $L/K$ as permutation group on the roots of $f$. Then conjugation of $G$ by an element of $N$ amounts to a renumbering of the roots of $f$. In particular, if $C_1, C_2$ are two conjugacy classes of $G$ fused in $N$, then whenever we have found an extension such that complex conjugation lies in class $C_1$, a simple renumbering provides an extension with complex conjugation in $C_2$. Thus in problem 2 we may restrict ourselves to consideration of classes of $G$ modulo the action of $N$. We have constructed extensions for all these possibilities up to degree 11 with the three above mentioned exceptions.

Problem 3 is completely solved up to degree 7. In degree 8 most transitive groups are covered but there are some primitive groups left where we cannot prove that we have found the minimal discriminant.

We close by giving a table containing some statistics about the number of polynomials in each degree. The # Classes column denotes the total number of classes of order 1 and 2 up to conjugation in the symmetric normalizer.

| Degree | # Groups | # Classes | # Polynomials |
|---:|---:|---:|---:|
| 2 | 1 | 2 | 549 |
| 3 | 2 | 3 | 619 |
| 4 | 5 | 13 | 2292 |
| 5 | 5 | 10 | 1489 |
| 6 | 16 | 48 | 5979 |
| 7 | 7 | 14 | 1360 |
| 8 | 50 | 233 | 14610 |
| 9 | 34 | 83 | 6144 |
| 10 | 45 | 184 | 12359 |
| 11 | 8 | 19 | 501 |
| 12 | 301 | 1895 | 43200 |
| 13 | 9 | 23 | 248 |
| 14 | 63 | 331 | 5155 |
| 15 | 104 | 395 | 4107 |

Universität Heidelberg, IWR, Im Neuenheimer Feld 368, 69120 Heidelberg, Germany.
*E-mail address*: `klueners@iwr.uni-heidelberg.de`

FB Mathematik/Informatik, Universität Kassel, Heinrich-Plett-Strasse 40, 34132 Kassel, Germany.
*E-mail address*: `malle@mathematik.uni-kassel.de`